\input amstex
\input amsppt.sty
\magnification=\magstep1
\hsize=32truecc
\vsize=22.5truecm
\baselineskip=16truept
\NoBlackBoxes
\TagsOnRight \pageno=1 \nologo
\def\Z{\Bbb Z}
\def\N{\Bbb N}

\def\l{\left}
\def\r{\right}
\def\bg{\bigg}
\def\({\bg(}
\def\[{\bg\lfloor}
\def\){\bg)}
\def\]{\bg\rfloor}
\def\t{\text}
\def\f{\frac}

\def\p{\ (\roman{mod}\ p)}

\def\bi{\binom}
\def\eq{\equiv}

\def\ls{\leqslant}
\def\gs{\geqslant}
\def\mo{\roman{mod}}

\def\ve{\varepsilon}

\def\M#1#2{\thickfracwithdelims[]\thickness0{#1}{#2}_6}

\def\m#1#2{\thickfracwithdelims\{\}\thickness0{#1}{#2}_m}

\def\Proof{\noindent{\it Proof}}

\def\Remark{\medskip\noindent{\it  Remark}}

\def\Ack{\medskip\noindent {\bf Acknowledgment}}
\hbox {J. Number Theory 160 (2016), 108--116.}
\bigskip
\topmatter
\title The least modulus for which consecutive polynomial values are distinct\endtitle
\author Zhi-Wei Sun\endauthor
\leftheadtext{Zhi-Wei Sun}
\rightheadtext{Zhi-Wei Sun}
\affil Department of Mathematics, Nanjing University\\
 Nanjing 210093, People's Republic of China
  \\  zwsun\@nju.edu.cn
  \\ {\tt http://math.nju.edu.cn/$\sim$zwsun}
\endaffil
\abstract Let $d\gs4$ and $c\in(-d,d)$ be relatively prime
integers. We show that for any
sufficiently large integer $n$ (in particular $n>24310$ suffices for $4\ls d\ls 36$),
the smallest prime $p\eq c\ (\mo\ d)$ with $p\gs(2dn-c)/(d-1)$ is the least positive integer $m$ with $2r(d)k(dk-c)\ (k=1,\ldots,n)$ pairwise distinct
modulo $m$, where $r(d)$ is the radical of $d$. We also conjecture that for any integer $n>4$
the least positive integer $m$ such that $|\{k(k-1)/2\ \mo\ m:\ k=1,\ldots,n\}|= |\{k(k-1)/2\ \mo\ m+2:\ k=1,\ldots,n\}|=n$
is the least prime $p\gs 2n-1$ with $p+2$ also prime.
\endabstract
\thanks 2010 {\it Mathematics Subject Classification}.\,Primary 11A41, 11B25;
Secondary 11A07, 11B75, 11N13, 11Y11.
\newline\indent Supported by the National Natural Science
Foundation (Grant No. 11571162) of China. The initial version was posted to arXiv as a preprint with the ID {\tt arXiv:1304.5988} in April 2013.
\endthanks
\endtopmatter
\document

\heading{1. Introduction}\endheading

To find nontrivial arithmetical functions taking only prime values is a fascinating topic in number theory.
In 1947 W. H. Mills [M]
showed that there exists a real number $A$ such that $\lfloor
A^{3^n}\rfloor$ is prime for every $n\in\Z^+=\{1,2,3,\ldots\}$;
unfortunately such a constant $A$ cannot be effectively found.

For each integer $h>1$ and sufficiently large integer $n$, it was determined in [BSW] the least positive integer $m$
with $1^h,2^h,\ldots,n^h$ pairwise distinct modulo $m$, but such integers $m$ are composite infinitely often.
In a recent paper [S] the author proved that
the smallest
integer $m>1$ such that $2k(k-1)$ mod $m$ for
$k=1,\ldots,n$ are pairwise distinct, is precisely the least prime greater than $2n-2$, and that
for $n\in\{4,5,\ldots\}$ the least positive integer $m$ such that $18k(3k-1)\ (k=1,\ldots,n)$ are pairwise distinct modulo $m$,
is the least prime $p>3n$ with $p\eq1\ (\mo\ 3)$.
When $d\in\{4,5,6,\ldots\}$ and $c\in(-d,d)$ are relatively prime, it is natural to ask whether there is a similar result
for primes in the arithmetic progression $\{c,c+d,c+2d,\ldots\}$ since there are infinitely many such primes by Dirichlet's theorem.

Based on our computation we discover the following general result.

\proclaim{Theorem 1.1} Let $d\gs4$ and $c\in(-d,d)$ be relatively prime
integers. Let
$$f_{d,c}(x):=2r(d)x(dx-c),\tag1.1$$
where $r(d)$ is the radical of $d$
$($i.e., the product of all the distinct prime divisors of $d)$.
For $n\in\Z^+$ define $m_{d,c}(n)$ as the least positive integer $m$ for which
the integers $f_{d,c}(k)\ (k=1,\ldots,n)$ are pairwise distinct modulo $m$.

{\rm (i)} If $n\in\Z^+$ is sufficiently large, then $m_{d,c}(n)$ is the least prime  $p\eq c\ (\mo\ d)$ with $p\gs(2dn-c)/(d-1)$.

{\rm (ii)} When $4\ls d\ls 36$ and $n>M_d$, the required result in the first part holds, where
$$\gather M_4=8,\ M_5=14,\ M_6=9,\ M_7=100,\ M_8=21,\ M_9=315,\ M_{10}=53,
\\M_{11}=1067,\ M_{12}=27,\ M_{13}=1074,\ M_{14}=122,\ M_{15}=809,\ M_{16}=329,
\\M_{17}=5115,\ M_{18}=95,\ M_{19}=5390,\ M_{20}=755,\ M_{21}=3672,\ M_{22}=640,
\\M_{23}=11193,\ M_{24}=220,\ M_{25}=12810,\ M_{26}=1207,\ M_{27}=7087,
\\M_{28}=2036,\ M_{29}=13250,\ M_{30}=177,\ M_{31}=24310,\ M_{32}=3678,
\\M_{33}=12794,\ M_{34}=5303,\ M_{35}=15628,\ M_{36}=551.
\endgather$$
\endproclaim
\Remark\ 1.1. To obtain the effective lower bounds $M_d\ (4\ls d\ls 36)$ in part (ii) of Theorem 1.1, we actually employ
some computational results of O. Ramar\'e and R. Rumely [RR] on primes in arithmetic progressions. Define
$$\gather c_4=-3,\ c_5=-1,\ c_6=1,\ c_7=-5,\ c_8=1,\ c_9=2,\ c_{10}=3,
\\c_{11}=-7,\ c_{12}=5,\ c_{13}=-5,\ c_{14}=-5,\ c_{15}=-1,\ c_{16}=11,
\\c_{17}=15,\ c_{18}=1,\ c_{19}=6,\ c_{20}=-9,\ c_{21}=1,\ c_{22}=5,
\\c_{23}=21,\ c_{24}=1,\ c_{25}=19,\ c_{26}=-3,\ c_{27}=23,
\\c_{28}=-9,\ c_{29}=-1,\ c_{30}=17,\ c_{31}=3,\ c_{32}=-1,
\\c_{33}=-5,\ c_{34}=15,\ c_{35}=12,\ c_{36}=23.
\endgather$$
Then, for every $d=4,\ldots,36$, the number $m_{d,c_d}(M_d)$ is {\it not} the least prime $p\eq c_d\ (\mo\ d)$ with $p\gs(2dM_d-c_d)/(d-1)$.
\medskip

 Theorem 1.1 with $d=4,5$ yields the following concrete consequence.

\proclaim{Corollary 1.1} {\rm (i)} For each integer $n\gs 6$, the least positive integer $m$ such that
$4k(4k-1)$ $($or $4k(4k+1))$ for $k=1,\ldots,n$ are pairwise distinct modulo $m$, is the least prime
$p\eq 1\ (\mo\ 4)$ with $p\gs(8n-1)/3$ $($resp., $p\eq-1\ (\mo\ 4)$ with $p\gs(8n+1)/3)$.

{\rm (ii)} Let $C_1=8,\ C_2=10,\ C_{-1}=15$ and $C_{-2}=5$.
For any $r\in\{\pm1,\pm2\}$ and integer $n\gs C_r$, the least positive integer $m$ such that
$10k(5k-r)$ for $k=1,\ldots,n$ are pairwise distinct modulo $m$, is the least prime
$p\eq r\ (\mo\ 5)$ with $p\gs (10n-r)/4$.
\endproclaim

As a supplement to Theorem 1.1, we are able to prove the following result for the cases $d=2,3$.

\proclaim{Theorem 1.2} For $d\in\{2,3\}$ and integer $c\in(-d,d)$, let $S_{d,c}$ be the set of all primes
$p\eq c\pmod d$ and powers of $d$. Then
$$\align m_{2,1}(n)=&\min\{a\gs 4n-1:\ a\in S_{2,1}\}\ \ \t{for}\ n\gs 5,
\\m_{2,-1}(n)=&\min\{a\gs 4n:\ a\in S_{2,-1}\}\ \ \t{for}\ n\gs 7,
\\m_{3,1}(n)=&\min\{a\gs3n:\ a\in S_{3,1}\}\  \ \t{for}\ n\gs 4,
\\m_{3,-1}(n)=&\min\{a\gs 3n:\ a\in S_{3,-1}\}\  \ \t{for}\ n\gs 5,
\\m_{3,2}(n)=&\min\{a\gs 3n-1:\ a\in S_{3,2}\}\  \ \t{for}\ n\gs 3,
\\m_{3,-2}(n)=&\min\{a\gs3n:\ a\in S_{3,-2}\}\  \ \t{for}\ n\gs 8.
\endalign$$
\endproclaim
\Remark\ 1.2. As the proof of Theorem 1.2 is quite similar to and even easier than that of Theorem 1.1, we omit the details of the proof.
Note that if $n>1$ is a power of two with $4n-1$ composite then $\min\{a\gs 4n-1:\ a\in S_{2,1}\}=4n$ is a power of two.
Also, if $n>1$ is a power of three then $\min\{a\gs 3n-1:\ a\in S_{3,2}\}=3n$ is a power of three.

\medskip

To conclude this section, we pose some new conjectures.

\proclaim{Conjecture 1.1} For any $d\in\Z^+$ there is a positive integer $n_d$ such that for any integer $n\gs n_d$
the least positive integer $m$ satisfying
$$\l|\l\{\bi k2\ \mo\ m:\ k=1,\ldots,n\r\}\r|=\l|\l\{\bi k2\ \mo\ m+2d:\ k=1,\ldots,n\r\}\r|=n$$
is the smallest prime $p\gs 2n-1$ with $p+2d$ also prime. Moreover, we may take
$$\gather n_1=5,\ n_2=n_3=6,\ n_4=10,\ n_5=9,
\\ n_6=8,\ n_7=9,\ n_8=18,\ n_9=11,\ n_{10}=9.\endgather$$
\endproclaim
\Remark\ 1.3. A well-known conjecture of de Polignac [P] asserts that
for any positive integer $d$ there are infinitely many prime pairs $\{p,q\}$ with $p-q=2d$.
\medskip

\proclaim{Conjecture 1.2} Let $n$ be any positive integer and consider the least positive integer $m$ such that
$$\l|\l\{\bi k2\ \mo\ m:\ k=1,\ldots,n\r\}\r|=\l|\l\{\bi k2\ \mo\ m+1:\ k=1,\ldots,n\r\}\r|=n.$$
Then, each of $m$ and $m+1$ is either a power of two $($including $2^0=1)$ or a prime times a power of two.
\endproclaim

\proclaim{Conjecture 1.3} Let $n$ be any positive integer.
Then the least positive integer $m$ of the form $x^2+x+1$ $($or $4x^2+1)$ with $x\in\Z$
such that the coefficients $\bi k2\ (k=1,\ldots,n)$ are pairwise distinct modulo $m$,
is the the smallest prime $p\gs 2n-1$ of the form  $x^2+x+1$ $($resp., $4x^2+1)$ with $x\in\Z$.
\endproclaim
\Remark\ 1.4. The conjecture that there are infinitely many primes of the form $x^2+x+1$ (or $4x^2+1$) is still open.
We may also replace $\bi k2$ in Conjecture 1.3 by $k(k-1)$.
\medskip

\proclaim{Conjecture 1.4} For any integer $n>2$, the smallest positive integer $m$ such that the integers $6p_k(p_k-1)\ (k=1,\ldots,n)$
are pairwise incongruent modulo $m$ is precisely the least prime $p\gs p_n$ dividing none of the numbers $p_i+p_j-1\ (1\ls i<j\ls n)$, where
$p_k$ denotes the $k$-th prime.
\endproclaim
\Remark\ 1.5. For any prime $p\gs p_n$ dividing none of the numbers $p_i+p_j-1\ (1\ls i<j\ls n)$,
clearly $p_j(p_j-1)-p_i(p_i-1)=(p_j-p_i)(p_i+p_j-1)\not\eq0\ (\mo\ p)$ for all $1\ls i<j\ls n$.

\medskip

We also have some other conjectures similar to Conjectures 1.1--1.4.

In the next section we provide some lemmas. Section 3 is devoted to our proof of Theorem 1.1.

\heading{2. Some lemmas}\endheading

\proclaim{Lemma 2.1} Let $c$ and $d>1$ be relatively prime integers. For any $\ve>0$, if $n\in\Z^+$ is large enough, then
there is a prime $p\eq c\ (\mo\ d)$ with
$$\f{d(2n-1)-c}{d-1}<p\ls\f{d((2+\ve)n-1)-c}{d-1}.$$
\endproclaim
\Proof. This is an easy consequence of the Prime Number Theorem for arithmetic progressions (cf. (1.5) of [CP, p.\,13] or Theorem 4.4.4 of [J, p.\,175])
which states that
$$|\{p\ls x:\ p\ \t{is a prime with}\ p\eq c\ (\mo\ d)\}|\sim\f x{\varphi(d)\log x}$$
as $x\to+\infty$, where $\varphi$ is Euler's totient function.  \qed

In view of (1.1), for any $c\in\Z$ and $d\in\Z^+$ we have the useful identity
$$f_{d,c}(l)-f_{d,c}(k)=2r(d)(l-k)(d(k+l)-c).\tag2.1$$

\proclaim{Lemma 2.2} Let $d>2$ and $c\in(-d,d)$ be relatively prime integers.
Suppose that  $p$ is a prime not exceeding $(d((2+\ve)n-1)-c)/(d-1)$ where $n\gs3d$ and $0<\ve\ls 2/(d-2)$.
Then
$$\aligned &f_{d,c}(k) \ (k=1,\ldots,n)\ \t{are pairwise distinct modulo}\ p
\\\iff&p\eq c\ (\mo\ d)\ \t{and}\ p>(d(2n-1)-c)/(d-1).
\endaligned\tag2.2$$
\endproclaim
\Proof. If $p\mid 2d$, then $p\mid 2r(d)$ and hence $f_{d,c}(k)\eq0\pmod p$ for all $k=1,\ldots,n$.
Note that $(d(2n-1)-c)/(d-1)\gs(3d-c)/(d-1)\gs2d/(d-1)>2$. If $p\mid d$ then $p\not\eq c\ (\mo\ d)$.
So (2.2) holds in the case $p\mid 2d$.

From now on we assume that $p\nmid 2d$. Then $jp\eq-c\ (\mo\ d)$ for some $1\ls j\ls d-1$.

Negating the right-hand side of (2.2), we suppose first that
$p\not\eq c\ (\mo\ d)$ or $p\ls(d(2n-1)-c)/(d-1)$. Write $jp+c=dq$ with $q\in\Z$. If $p\not\eq c\ (\mo\ d)$, then $j\ls d-2$ and hence
$$\align q\ls&\f cd+\f{d-2}dp\ls\f cd+\f{d-2}d\cdot\f{d((2+\ve)n-1)-c}{d-1}
\\\ls&\f {c-d(d-2)}{d(d-1)}+\f{d-2}{d-1}\l(2+\f2{d-2}\r)n<2n.
\endalign$$
If $p\eq c\ (\mo\ d)$ and $p\ls(d(2n-1)-c)/(d-1)$, then $j=d-1$ and $q\ls 2n-1$.
When $q>2$, we have $0<k:=\lfloor(q-1)/2\rfloor<l:=\lfloor(q+2)/2\rfloor\ls n$, also
$$d(k+l)-c=dq-c=jp\eq0\ \ (\mo\ p)$$
and hence $f_{d,c}(k)\eq f_{d,c}(l)\pmod p$ in view of (2.1).
If $q\ls2$, then $p\ls jp=dq-c\ls2d-c<3d\ls n$ and $f_{d,c}(p+1)\eq f_{d,c}(1)\pmod p$.

We now assume the right-hand side of (2.2). Then $(d-1)p+c=dq$ for some integer $q\gs2n$.
In view of (2.1), we only need to show that $p\nmid (l-k)$ and $d(k+l)\not\eq c\pmod p$ for any $1\ls k<l\ls n$.
Note that
$$0<l-k<n\ls\f{dq}{2d}=\f{(d-1)p+c}{2d}<\f{p+1}2\ls p$$
and also $d(k+l)-c\ls d(2n-1)-c<(d-1)p$. If $d(k+l)\eq c\pmod p$, then for some $t=1,\ldots,d-2$ we have $d(k+l)-c=tp\eq tc\not\eq-c\ (\mo\ d)$,
which leads to a contradiction.

The proof of Lemma 2.2 is now complete. \qed

\proclaim{Lemma 2.3} Let $d>2$ and $c\in(-d,d)$ be relatively prime integers, and let $n\gs6d$ be an integer.
Suppose that $m\in[n,(d((2+\ve)n-1)-c)/(d-1)]$ is a power of two or twice an odd prime,
where $0<\ve\ls 2/3$. Then, there are $1\ls k<l\ls n$ such that
$f_{d,c}(k)\eq f_{d,c}(l)\pmod m$.
\endproclaim
\Proof. Note that $m\gs n\gs 6d>4$ and
$$\f m4\ls\f{d((2+\ve)n-1)-c}{4(d-1)}<\f{d(2+\ve)}{4(d-1)}n\ls\f{d(2+2/3)}{4(d-1)}n=\f{8dn}{8d+4(d-3)}\ls n.$$

If $d$ is even and $m$ is a power of two, then for $k=1$ and $l=m/4+1\ls n$ we have $m\mid 2r(d)(l-k)$ and hence
$f_{d,c}(k)\eq f_{d,c}(l)\pmod m$ by (2.1).
If $m=2p$ with $p$ an odd prime dividing $d$, then $m\mid 2r(d)$ and hence $f_{d,c}(k)\eq0\pmod m$
for all $k=1,\ldots,n$.

In the other cases, $d$ and $m/2$ are relatively prime. Thus $jd\eq c\ (\mo\ m/2)$ for some $j=1,\ldots,m/2$.
If $j\ls 2$, then
$$\f m2\ls jd-c\ls 2d-c<3d\ls\f n2$$
which contradicts $m\gs n$. So $3\ls j\ls m/2$ and hence
$$0<k:=\l\lfloor\f{j-1}2\r\rfloor<l:=\l\lfloor\f{j+2}2\r\rfloor\ls\f{m}4+1<n+1.$$
Since $d(k+l)-c=jd-c\eq0\ (\mo\ m/2)$, by (2.1) we have $f_{d,c}(k)\eq f_{d,c}(l)\pmod m$.
This concludes the proof. \qed

\heading{3. Proof of Theorem 1.1}\endheading

\medskip\noindent{\it Proof of Theorem 1.1}.  Let $\ve=2/(\max\{11,d\}-2)$. By Lemma 2.1,
if $n\in\Z^+$ is large enough then there is at least a prime $p\eq c\ (\mo\ d)$ with
$$\f{d(2n-1)-c}{d-1}<p\ls\f{d((2+\ve)n-1)-c}{d-1}.\tag3.1$$

(i) Choose an integer $N\gs \max\{6d,243\}$ such that for any integer $n\gs N$  there is a prime $p\eq c\pmod d$ satisfying (3.1).
Fix an integer $n\gs N$ and let $m=m_{d,c}(n)$.
Clearly $m\gs n$. By Lemma 2.2, $m\ls m'$ where $m'$ denotes the least prime $p\eq c\ (\mo\ d)$ satisfying (3.1).

Assume that $m\not=m'$. We want to reach a contradiction.
Clearly $m$ is not a prime by Lemma 2.2. Note that  $\ve\ls 2/9$. In view of Lemma 2.3,
$m$ is neither a power of two nor twice an odd prime.
So we have $m=pq$ for some odd prime $p$ and integer $q>2$.
Observe that
$$\f m3\ls\f{d((2+\ve)n-1)-c}{3(d-1)}
<\f {d(2+2/9)}{3(d-1)}n=\f{20d}{27(d-1)}n\ls \f{80}{81}n$$
and hence
$$\f m3+3<\f{80}{81}n+\f n{81}=n.\tag3.2$$

If $p\mid d$, then for $k:=1$ and $l:=q+1=m/p+1<m/3+3<n$, we have $pq\mid r(d)(l-k)$ and hence
$f_{d,c}(k)\eq f_{d,c}(l)\pmod m$ by (2.1).

Now suppose that $p\nmid d$. Then $2dk\eq c-dq\ (\mo\ p)$ for some $1\ls k\ls p$. Clearly, $l:=k+q\ls p+q=m/q+m/p$.
Note that
$$(l-k)(d(l+k)-c)=q(d(2k+q)-c)\eq 0\pmod{pq}$$
and hence $f_{d,c}(k)\eq f_{d,c}(l)\pmod m$ by (2.1).
If $\min\{p,q\}\ls 4$, then
$$l\ls p+q=\f m{\min\{p,q\}}+\min\{p,q\}\ls\f m3+4<n+1$$
by (3.2). If $\min\{p,q\}\gs5$, then
$$l\ls\f mq+\f mp\ls \max\l\{\f m6+\f m7,\ \f m5+\f m8\r\}<\f m3<n$$
since $pq=m\gs n\gs 243\gs40$.
So we get a contradiction as desired.

(ii) Now assume that $4\ls d\ls 36$.
By Table 1 of [RR, p.\,419], we have
$$(1-\ve_d)\f x{\varphi(d)}\ls \theta(x;c,d)\ls (1+\ve_d)\f x{\varphi(d)}\ \ \ \ \t{for all}\ x\gs 10^{10},\tag3.3$$
where
$$\theta(x;c,d):=\sum\Sb p\ls x\\p\eq c\pmod d\endSb\log p\quad \ \t{with}\ p\ \t{prime},$$
and
$$\align \ve_4=&0.002238,\ \ve_5=0.002785,\ \ve_6=0.002238,\ \ve_7=0.003248,\ \ve_8=0.002811,
\\\ve_9=&0.003228,\ \ve_{10}=0.002785,\ \ve_{11}=0.004125,\ \ve_{12}=0.002781,\ \ve_{13}=0.004560,
\\\ve_{14}=&0.003248,\ \ve_{15}=0.008634,\ \ve_{16}=0.008994,\ \ve_{17}=0.010746,\ \ve_{18}=0.003228,
\\\ve_{19}=&0.011892,\ \ve_{20}=0.008501,\ \ve_{21}=0.009708,\ \ve_{22}=0.004125,\ \ve_{23}=0.012682,
\\\ve_{24}=&0.008173,\ \ve_{25}=0.012214,\ \ve_{26}=0.004560,\ \ve_{27}=0.011579,\ \ve_{28}=0.009908,
\\\ve_{29}=&0.014102,\ \ve_{30}=0.008634,\ \ve_{31}=0.014535,\ \ve_{32}=0.011103,\ \ve_{33}=0.011685,
\\\ve_{34}=&0.010746,\ \ve_{35}=0.012809,\ \ve_{36}=0.009544.
\endalign$$
As $\ve=2/(\max\{11,d\}-2)$, we can easily verify that
$$\f{\ve}2-\f 2{10^{10}}>\f{2\ve_d}{1-\ve_d}=\f{1+\ve_d}{1-\ve_d}-1.$$
If $n\gs 10^{10}/2$, then
$$((2+\ve)n-2)\f d{d-1}\gs 2n\f d{d-1}>10^{10}$$
and $$\f{\ve}2-\f1n+1\gs\f{\ve}2-\f 2{10^{10}}+1>\f{1+\ve_d}{1-\ve_d},$$
hence by (3.3) we have
$$\align&\f{\theta(((2+\ve)n-2)d/(d-1);c,d)}{\theta(2nd/(d-1);c,d)}
\\\gs&\f{(1-\ve_d)((2+\ve)n-2)d/(d-1)}{(1+\ve_d)2nd/(d-1)}
=\f{1-\ve_d}{1+\ve_d}\l(1+\f{\ve} 2-\f1n\r)>1
\endalign$$
and thus there is a prime $p\eq c\ (\mo\ d)$ for which
$$\f{2dn}{d-1}<p\ls\f{((2+\ve)n-2)d}{d-1}$$
and hence (3.1) holds.

Let $N_d$ be the least positive integer such that for any $n=N_d,\ldots,10^{10}/2$ and any $a\in\Z$ relatively prime to $d$, the interval
$(2dn/(d-1),((2+\ve)n-2)d/(d-1))$ contains a prime congruent to $a$ modulo $d$. Via a computer we find that
$$\gather N_4=79,\ N_5=206,\ N_6=103,\ N_7=471,\ N_8=301,\ N_9=356, N_{10}=232,
\\N_{11}=1079,\ N_{12}=346,\ N_{13}=1166,\ N_{14}=806,\ N_{15}=1310,\ N_{16}=2183,
\\N_{17}=5153,\ N_{18}=1135,\ N_{19}=5402,\ N_{20}=2388,\ N_{21}=4059,\ N_{22}=2934,
\\N_{23}=11246,\ N_{24}=2480,\ N_{25}=13144,\ N_{26}=4775,\ N_{27}=11646,
\\\ N_{28}=5314,\ N_{29}=13478,\ N_{30}=5215,\ N_{31}=24334,\ N_{32}=8964,
\\\ N_{33}=15044,\ N_{34}=14748,\ N_{35}=16896,\ N_{36}=9847.
\endgather$$
For any integer $n\gs N_d$, there is a prime $p\eq c\pmod d$ satisfying (3.1).
Note that $6d\ls 6\times 36<243$.
So, for $n\gs N=\max\{N_d,243\}$ we may apply part (i) to get the desired result.
If $M_d<n\ls\max\{N_d,243\}$, then we can easily verify the desired result via a computer.

In view of the above, we have completed the proof of Theorem 1.1. \qed

\Ack. The author would like to thank the referee for helpful comments.

\enddocument